\documentclass[a4paper,11pt]{article}
\usepackage{amsfonts}
\usepackage{amsmath}
\usepackage{amssymb}

\usepackage[left=1in,top=1in,right=1in,bottom=1in]{geometry}

\newcommand{\dirac}{I}               

\usepackage{setspace}

\newcommand{\EE}{\mathbb{E}}                    
\newcommand{\II}{\mathbb{I}}                    
\newcommand{\steinoperator}{{\mathcal{A}}}      
\newcommand{\distributionof}[1]{{\mathcal{L}(#1)}} 
\newcommand{\sivn}{\sum_{i=1}^n}               
\newcommand{\skn}{\sum_{k=1}^n}                 

\newcommand{\cdik}[3]{\big\{\Delta^{#1}g(#2)\vert\II_i=#3\big\}}
\newcommand{\isdefinedas}{:=}                   

\newcommand{\ewai}{\EE\ab{\widetilde W^{(i)}-W^{(i)}}}
\newcommand{\deltdui}[2]{
\big\{ \Delta g(W^{(i)}+#1)\vert\II_i=#2\big\}}

\newtheorem{theorem}{\bf Theorem}[section]
\newtheorem{proposition}[theorem]{\bf Proposition}
\newtheorem{corollary}[theorem]{\bf Corollary}
\newtheorem{remark}[theorem]{\bf Remark}
\newtheorem{lemma}[theorem]{\bf Lemma}

\newcommand{\qed}{\hfill \ensuremath{\Box}}

\newcommand{\norm}[1]{\|#1\|}                    
\newcommand{\Norm}[1]{\Big\|#1\Big\|}            
\newcommand{\ab}[1]{\vert#1\vert}                
\newcommand{\Ab}[1]{\Big\vert#1\Big\vert}        

\newcommand{\ee}{{\mathrm e}}

\newcommand{\exponent}[1]{\exp\{#1\}}            
\newcommand{\Exponent}[1]{\exp\Bigl\{#1\Bigr\}}  

\newcommand{\Expect}{\mathbb{E}}                 

\newcommand{\asta}{{\cal A}}                      

\newcommand{\eL}{{\cal L}}
\newcommand{\F}{{\cal F}}

\newcommand{\Zplus}{{\mathbb Z}_+}
\newcommand{\RR}{\mathbb{R}}
\newcommand{\NN}{\mathbb{N}}

\newcommand{\binomial}[2]{\genfrac{(}{)}{0pt}{}{#1}{#2}}

\newcommand{\Proof}{\textbf{Proof. }}            
\newcommand{\Proofof}[1]{\textbf{Proof of #1. }} 

\newcommand{\NB}{\mathrm{NB}}           
\newcommand{\BIP}{\mathrm{BCP}}    

\newcommand{\qbar}{\bar{q}}
\newcommand{\pbar}{\bar{p}}
\title{On Stein Operators for Discrete  Approximations}

\author{ N. S. Upadhye,  V. \v Cekanavi\v cius and P. Vellaisamy\\
{\small Department of Mathematics, Indian Institute of Technology
Madras,}\\
{\small Chennai - 600036, India }\\{\small  E-mail:
neelesh@iitm.ac.in}
\\{\small
Department of Mathematics and Informatics, Vilnius University,}\\
{\small Naugarduko 24, Vilnius - 03225, Lithuania,}\\{\small
E-mail: vydas.cekanavicius@mif.vu.lt}\\{\small and} {\small
 Department of Mathematics, Indian Institute of Technology Bombay,} \\
 {\small Powai, Mumbai -
400076, India.}\\{\small  E-mail: pv@math.iitb.ac.in} }
\date{}
\begin{document}

\maketitle

\begin{abstract}

\noindent
In this paper, a new method based on probability generating functions is used to obtain multiple Stein operators for various random variables closely related to Poisson, binomial and negative binomial distributions. Also, Stein operators for certain compound distributions, where the random summand satisfies Panjer's recurrence relation, are derived. A well-known perturbation approach for Stein's method is used to obtain total variation bounds for the distributions mentioned above. The importance of such approximations is illustrated, for example, by the binomial convoluted with Poisson approximation to sums of independent and dependent indicator random variables.

\vspace*{.5cm} \noindent {\emph{Keywords:} \small  Binomial
distribution, compound Poisson distribution, perturbation, Panjer's recursion, Stein's method,
total variation norm.}

\vspace*{.5cm} \noindent {\small {\it MSC 2000 Subject
Classification}:
Primary 60F05.   
Secondary  62E20;     
}

\end{abstract}

\section{Introduction}
\label{Intro}
Stein's method is known to be one of the powerful techniques for
probability approximations and there is a vast literature available on this topic. For  details and applications of Stein's method,
 see Barbour {\it et al.~}(1992b), Chen {\it et al.~}(2011), Ross
(2011) and Norudin and Peccati (2012).
For some recent developments, see
Eichelsbacher and Reinert (2008), Chen {\it et al.~}(2011), Norudin and
Peccati (2012), Daly, Lefevre and Utev (2012), Ley and Swan (2013a, 2013b)
and the references therein. The
method is based on the construction of a characteristic operator for an
approximation problem. Different approaches are used for deriving Stein operators (see, Reinert (2005)). For instance, a Stein operator can be treated in
the framework of birth-death processes (Brown and Xia (2001)).
Stein's method for discrete distributions has been  independently and simultaneously developed by  Goldstein and Reinert (2013) and Ley and Swan (2013). More recently, Ley {\it et al.~}(2014) has proposed a canonical operator, for both continuous and discrete distributions, and a general approach to obtain bounds on approximation problems.

In this paper, we consider the random variables (rvs) concentrated on $\Zplus=\{0,1,2,\dots\}$ with distributions having the form of convoluted measures or random sums. Using their probability
generating functions ($PGF$'s), we derive Stein operators for discrete probability approximations.
In particular, the existence of multiple Stein operators (in the case of convoluted measures)
for an approximation problem is shown and the corresponding bounds are derived, using perturbation technique, and compared for the case of indicator rvs. Although the existence of infinite families of Stein operators for many common distributions is already well-known (see, Goldstein and Reinert (2005) and Ley {\it et al.~}(2014)), this comparison may benefit the readers, as it is illustrated for the first time (in case of convoluted measures) to the best of our knowledge.

Next, we describe a typical procedure for Stein's method on $\Zplus$-valued rvs.
Let $Y$ be a $\Zplus$-valued rv with $\Expect (|Y|) < \infty$, ${\cal F} := \{f| f:\Zplus \to \RR {\mbox{~and  is bounded}} \}$ and ${\cal G}_Y = \{g \in {\cal F}| g(0)=0~{\rm and}~ g(x) = 0~{\rm for}~ x \notin supp(Y)\}$, where $supp(Y)$ denotes the support of rv $Y$.
We want to bound $\Expect f(Z)-\Expect f(Y)$ for
some rv $Z$ concentrated on $\Zplus$ and $f \in {\cal F}$. Stein's method is then realized in three
consecutive steps. First, for any
$g \in {\cal G}_Y$,  a linear operator $\asta$ satisfying
$\Expect(\asta g)(Y)=0$ is established and is called a
Stein operator. For a general framework of Stein operators, the reader is referred to Stein (1986), Stein {\it et al.~}(2004), D\"{o}bler (2012), Goldstein and Reinert (2013), Fulman and Goldstein (2014), Ley and Swan (2013a, 2013b) and Ley {\it et al.~}(2014).

\noindent
In the next step, the so-called Stein equation
\begin{equation}
 (\asta g)(j) = f(j)-\Expect f(Y), \quad j \in
\Zplus, f \in {\cal F}\label{solo}
\end{equation}
is solved with respect to $g(j)$ in terms of $f$ and is referred to as a solution to the Stein equation (\ref{solo}). As a rule, solutions to the
Stein equations have useful properties, such as $\norm{\Delta
g}:=\sup_{j\in\Zplus}\ab{\Delta g(j)}$ is small, where $\Delta
g(j):=g(j+1)-g(j)$ denotes the first forward difference. Note that the properties of $\Delta g$ depend on the form of $\asta$ and some properties of $Y$. Finally, taking expectations on both sides of (\ref{solo}), we get
\begin{equation}
 \Expect f(Z)-\Expect f(Y)= \Expect (\asta g)(Z)\label{se}
\end{equation}
and bounds for $\Expect(\asta g)(Z)$ are established through the bounds for $\Delta g $ and
$\Delta^{k+1}g(j):=\Delta^{k}(g(j+1)-g(j))$, $k =1,2,\ldots$ .
For more details on Stein's method under a general setup, we refer the readers to Goldstein and Reinert (2005, 2013), Ley {\it et al.~} (2014), Barbour and Chen (2014) and the references therein.


\noindent
For some standard distributions, a Stein operator can be established easily. Indeed, let
$\mu_j:=P(Y=j)>0,  j \in \Zplus$. Then
$\sum_{j=0}^\infty
\mu_j\bigg(\frac{(j+1)\mu_{j+1}}{\mu_j}g(j+1)-jg(j)\bigg)=0.$
Therefore,
 \begin{equation}
(\asta g)(j)= \frac{(j+1)\mu_{j+1}}{\mu_j}g(j+1)-jg(j), \quad
j\in\Zplus, \label{asol}
\end{equation}
and it can be easily verified that ${\Expect (\asta g) (Y)} = 0$.
Some well-known examples are listed below.
\begin{itemize}
\item[1)] For $\alpha >0$, let $Y_1$ be a Poisson $P(\alpha)$ rv with
$\mu_j= P(Y_1=j)=\alpha^j\ee^{-\alpha}/j!$. Then
\begin{equation}
(\asta g)(j)=\alpha g(j+1)-j g(j), \quad j\in\Zplus.
\label{astpois}
\end{equation}
\item[2)] Let $0<p<1$,  $q=1-p$, $\widetilde
M >1$, and $Y_2$ have the pseudo-binomial distribution (see \v
Cekanavi\v cius and Roos (2004), p.~370) so that
\[\mu_j= P(Y_2=j)=\frac{1}{\widetilde{C}}\binomial{\widetilde{M}}{j}
   p^jq^{\widetilde{M}-j}, \quad j\in\{0,1,\dots,
   \lfloor\widetilde{M}\rfloor\},\]
where
$\widetilde{C}=\sum_{j=0}^{\lfloor\widetilde{M}\rfloor}
   \binomial{\widetilde{M}}{j}p^jq^{\widetilde{M}-j},$
 $\lfloor\widetilde M\rfloor$ denotes integer part of $\widetilde M$
and
$\binomial{\widetilde{M}}{j}=\frac{\widetilde M(\widetilde M-1)\cdots (\widetilde
M-j+1)}{j!}.$
 If $\widetilde M$ is an integer, then $Y_2$ is a binomial rv. Suppose
now  $g(0)=0$ and $g(\lfloor\widetilde{M}\rfloor
+1)=g(\lfloor\widetilde{M}\rfloor+2)=...=0$. Then, from (\ref{asol})
\[
(\asta g)(j) = \frac{(\widetilde{M}-j)p}{q}g(j+1)-jg(j), \quad
j=0,1,\dots\lfloor \widetilde M\rfloor.
\]
\noindent Multiplying the above expression by $q$, we can get  the following
Stein operator:
\begin{equation}
(\asta g)(j) = (\widetilde{M}-j)pg(j+1)-jqg(j), \quad
j=0,1,\dots\lfloor \widetilde M\rfloor. \label{astbin}
\end{equation}

\item[3)] Let $Y_3 \sim \NB(r, \pbar)$, $0<\pbar<1$,  be negative
binomial distribution with
 $\mu_j = P(Y_3=j)=\Gamma(r+j)/(\Gamma(r)j!) \pbar^r \qbar^j$, for $j\in\Zplus$, $r>0$ and
 $\qbar=1-\pbar$.
Then (\ref{asol}) reduces to
\begin{equation}
(\asta g)(j):= \qbar (r +j) g(j+1) - j g(j), \quad j\in\Zplus.
\label{astnb}
\end{equation}
\end{itemize}

Observe that equation (\ref{asol}) is not that useful if we do not have
simple expressions for $\mu_j$ and especially for $\mu_{j}/\mu_{j+1}.$
One such class is the  Ord family of cumulative distributions
(see Afendras {\it et al.~}(2014)).
 For example, if we consider
compound distribution or convolution of two or more distributions,
then $\mu_j$'s are usually expressed through sums or converging
series of probabilities. Therefore, some other refined approaches
for obtaining Stein operator(s) are needed.

The  paper is organized as follows. In Section
\ref{Soperator},  we use $PGF$ approach to obtain general expressions for Stein operators arising out of convolution of rvs and random sums that satisfy Panjer's recursive relation.
These operators are then seen as perturbations of known operators for standard distributions which motivate the discussion about perturbation approach and its applications.
In Section \ref{Sperturbations}, some facts about the perturbation approach
to a solution of Stein equation are discussed and applied to operators derived in Section \ref{Soperator}. In Section \ref{AIndicator}, as an application,
an approximation problem for the distribution of the sum of possibly dependent indicator variables by the
convolution of Poisson and binomial distribution is considered. We show that such
approximations can be treated either as Poisson perturbation or as
binomial perturbation, leading to two different bounds.
Finally, we mention that though the approach is restricted to distributional
approximations, its ideas can be extended for
approximations to signed measures as well.

\section{Stein Operators via $PGF$}\label{Soperator}
In this section, the $PGF$ approach is used to derive operators satisfying  $\Expect(\asta g)(Y)=0$ for $g\in{\cal G}_Y$.
The construction of $\asta$ is well-known if probabilities
of approximating distribution satisfy some recursive relation and it can be easily verified by using this approach.  Indeed, the $PGF$ has been used as a tool for establishing  Panjer's recurrence relations; see, for example, Sundt (1992) and Hess {\it et al.~}(2002). Note also that, strictly speaking,  $\asta$ can be called a Stein operator only if it is used in (\ref{se}) with $g$ satisfying (\ref{solo}).  Moreover, one expects $g$ to have some useful properties.
In Section 3, we show that,  the majority of operators considered below have solutions to (1) with properties  typical for the Stein method.

Next, we demonstrate how the $PGF$ approach can be used to derive the Stein operators for compound Poisson distribution, certain convolution of distributions and a compound distribution where the summand  satisfy the Panjer's recurrence relation.

\subsection{The General Idea}
\noindent
Let $N$ be  a $\Zplus$-valued
rv with $\mu_k=P(N=k)$ and finite mean. Then its
$PGF$
\begin{equation}
G_N(z)=\sum_{k=0}^\infty\mu_kz^k \label{PGFN}
\end{equation}
satisfies
\begin{equation}
G_{N}^{'}(z)=\frac{d}{dz}G_N(z)=\sum_{k=1}^\infty
k\mu_kz^{k-1}=\sum_{k=0}^\infty (k+1)\mu_{k+1}z^k, \label{dmz}
\end{equation}
where prime denotes the derivative with respect to $z$. If we can express $G_N^{'}(z)$ through $G_N(z)$ then, by
collecting factors corresponding to  $z^k$, the recursion follows. One can easily verify the Stein operators derived for standard distributions in the previous section, using this approach.

Next, we demonstrate the $PGF$ approach to derive a Stein operator for compound Poisson
distribution.

Let $\{X_j\}$ be an iid sequence of random variables with $P(X_j = k) = p_k$, $k=0,1,2,\ldots$ . Also, let $N \sim P(\lambda)$ and  be independent of the $\{X_j\}$. Then the distribution of $Y_4 := \sum_{j=1}^N X_j$ is known as compound Poisson distribution with the $PGF$
\begin{equation}
G_{cp}(z)=\Exponent{\sum_{j=1}^\infty\lambda_j(z^j-1)},\label{CPPGF}
\end{equation}
where $\lambda_j = \lambda p_j$ and $\sum_{j=1}^\infty j \ab{\lambda_j}<\infty$.
Then
\[
G_{cp}^{'}(z)=G_{cp}(z)\sum_{j=1}^\infty j\lambda_j z^{j-1}=
\sum_{k=0}^\infty\mu_kz^k\sum_{j=1}^\infty j\lambda_j
z^{j-1}=\sum_{k=0}^\infty
z^k\sum_{m=0}^k\mu_m(k-m+1)\lambda_{k-m+1}.
\]
Comparing the last expression to the right-hand side of (\ref{dmz}), we
obtain the recursive relation, for all $k\in\Zplus$, as
\[\sum_{m=0}^k\mu_m(k-m+1)\lambda_{k-m+1}-(k+1)\mu_{k+1}=0.\]
Then, for $g \in {\cal G}_{Y_4}$, we have
\begin{eqnarray*}
0&=&\sum_{k=0}^\infty
g(k+1)\Big[\sum_{m=0}^k\mu_m(k-m+1)\lambda_{k-m+1}-(k+1)\mu_{k+1}\Big]\\
&=& \sum_{m=0}^\infty\mu_m\Big[\sum_{k=m}^\infty
g(k+1)(k-m+1)\lambda_{k-m+1}-m
g(m)\Big]\\
&=&\sum_{m=0}^\infty\mu_m\Big[\sum_{j=1}^\infty j\lambda_j g(j+m)-
m g(m)\Big].
\end{eqnarray*}
Therefore, a Stein operator for compound Poisson distribution, defined in (\ref{CPPGF}), is
\begin{eqnarray}
(\asta g)(j) &=&\sum_{l=1}^\infty
l\lambda_lg(j+l)-jg(j)\nonumber\\
&=&\sum_{l=1}^\infty l\lambda_l g(j+1)-jg(j)+\sum_{m=2}^\infty
m\lambda_m\sum_{l=1}^{m-1}\Delta g(j+l) ,\quad j\in\Zplus,
\label{opCP}
\end{eqnarray}
since ${\Expect (\asta g)(Y_4)} = 0$.
This operator coincides with the one from Barbour  {\it et al.~}(1992a).

\noindent
Next, we derive multiple Stein operators for convolution of standard distributions discussed above.

\subsection{Convolutions of Distributions}
Recall that  $Y_1\sim P(\alpha)$ ($\alpha>0$), $Y_2\sim  Bi(M,p)$
($M\in\NN$, $0<p<1$), $Y_3\sim \NB(r,\pbar)$ ($0<\pbar<1$,
$r>0$) and $Y_4$ follows the compound Poisson distribution defined in (\ref{CPPGF}). We assume that $Y_1,Y_2$,$Y_3$ and $Y_4$ are independent. Then the $PGF$'s of
$Y_1+Y_2$, $Y_2$ and $Y_3$ are given by
\begin{equation}
 G_{12}(z)=(q+pz)^M\exponent{\alpha (z-1)},\quad
G_2(z)=(q+pz)^M,\quad G_3(z)=\bigg(\frac{\pbar}{1-\qbar
z}\bigg)^r, \label{a3}
\end{equation}
respectively. Here $\qbar=1-\pbar$ and $q=1-p$. We now derive Stein operators for the convolutions of various combinations of $Y_1, Y_2, Y_3$ and $Y_4$.

\begin{proposition} \label{pro1} Let $G_{cp}(z)$ be the $PGF$ of $Y_4$ and  $\lambda=\sum_{j=1}^\infty j\lambda_j$. Then we have the following results:\\
 (i) The rv $Y_{24} = Y_2+Y_4$ has the PGF
$G_2(z)G_{cp}(z)$  and its Stein operator, for $g\in {\cal G}_{Y_{24}}$, is
\begin{eqnarray}
(\asta g)(j)&=&\bigg(M+\frac{\lambda}
{p}-j\bigg)pg(j+1)-qjg(j)\nonumber\\
&&
+\sum_{m=2}^\infty\big(qm\lambda_m+p(m-1)\lambda_{m-1}\big)\sum_{l=1}^{m-1}\Delta
g(j+l). \label{a1}
\end{eqnarray}
(ii) The rv $Y_{34}=Y_3 + Y_4$ has the PGF $G_3(z)G_{cp}(z)$  and has a Stein operator, for $g\in {\cal G}_{Y_{34}}$,
\begin{eqnarray}
(\asta g)(j)&=&\bigg(\frac{\lambda
\pbar}{\qbar}+r+j\bigg)\qbar g(j+1)-jg(j)\nonumber\\
&&+\sum_{m=2}^\infty\big(m\lambda_m-\qbar(m-1)\lambda_{m-1}\big)\sum_{l=1}^{m-1}\Delta
g(j+l). \label{a2}
\end{eqnarray}
\end{proposition}

\noindent \textbf{Proof.} Write $G_2(z)G_{cp}(z)=\sum_{k=0}^\infty
\mu_kz^k$. Differentiating with respect to $z$, we get the identity
\[
\sum_{k=0}^\infty\mu_kz^k\Big(\frac{Mp}{q+pz}+\sum_{j=1}^\infty\lambda_jjz^{j-1}\Big)=\sum_{k=0}^\infty
k\mu_kz^{k-1}.
\]
Multiplying both sides by $(q+pz)$ and collecting the terms corresponding to $z^k$,
we obtain the recursive relation
\[
\sum_{m=0}^k\mu_m(q\lambda_{k-m+1}(k-m+1)+p(k-m)\lambda_{k-m})-(k+1)\mu_{k+1}q+(Mp-pk)\mu_k=0.
\]
Multiplying the last equation by $g(k+1)$ and summing over all
nonnegative  integer $k$ leads to (\ref{a1}).

\noindent
To prove (\ref{a2}),  let
$G_3(z)G_{cp}(z)=\sum_{k=0}^\infty \mu_kz^k$. Differentiating with
respect to $z$ gives the identity
\[
\sum_{k=0}^\infty\mu_kz^k\Big(\frac{r\qbar}{1-\qbar
z}+\sum_{j=1}^\infty\lambda_jjz^{j-1}\Big)=\sum_{k=0}^\infty
k\mu_kz^{k-1}.
\]
Multiplying both sides by $(1-\qbar z)$ and collecting the terms corresponding to $z^k$,
we obtain
\[
\sum_{m=0}^k\mu_m(\lambda_{k-m+1}(k-m+1)-\qbar(k-m)\lambda_{k-m})-(k+1)\mu_{k+1}+\qbar(k+r)\mu_k=0.
\]
Multiply the above equation by $g(k+1)$ and then sum
over $k\in \Zplus $ to obtain the result. \qed

\begin{proposition} \label{binp} Let $Y_{12}=Y_1+Y_2$ have PGF $G_{12}(z)$ as defined in
(\ref{a3}).
 Then, for $j\in\Zplus$ and
$g \in \mathcal{G}_{Y_{12}}$, a Stein operator for $Y_{12}$ is
\begin{equation}
(\asta\, g)(j)=(Mp+\alpha -j p)g(j+1)-jqg(j)+p\alpha\Delta g(j+1).
\label{astaB}
\end{equation}
If in addition $p<q$, then
\begin{equation}
(\asta g)(j)=(\alpha+Mp)g(j+1)-jg(j)+M\sum_{l=2}^\infty
(-1)^{l+1}\bigg(\frac{p}{q}\bigg)^l\sum_{k=1}^{l-1}\Delta g(j+k).
 \label{astaCP}
\end{equation}
\end{proposition}

\noindent \Proof Observe that (\ref{astaCP}) follows from  (\ref{opCP}) and the expansion
\begin{equation}
(q+pz)^M=\Exponent{M\sum_{i=1}^\infty\frac{(-1)^{i+1}}{i}\bigg(\frac{p}{q}\bigg)^i(z^i-1)}.
\label{uu1}\end{equation} Note that (\ref{astaB}) is a special case
of (\ref{a1}).\qed

\begin{remark}
(i) As is known in the literature (see Goldstein and Reinert (2005)), we
 have two significantly different Stein operators (see (\ref{astaB}) and (\ref{astaCP})) for the approximation
problem.

\noindent (ii)Observe that, the operator given in (\ref{astaB}) is similar to the operator given in (\ref{astbin}), where $\widetilde{M}$ is replaced by $M+ \alpha/p$, except for the last term, and hence is  known as a binomial perturbation.

\noindent (iii) Similarly, the operator given in (\ref{astaCP}) is similar to the operator given in (\ref{astpois}), where $\alpha$ is replaced by $Mp + \alpha$, except for the last sum, leading to a Poisson perturbation.
\end{remark}

\noindent
Next, we demonstrate
that the number of such operators might be even larger.  We
consider the convolution of negative binomial and binomial
distributions. It is logical to use the binomial approximation for
sums  of rv's with variances smaller than their means and the negative
binomial approximation if variances are larger than means.
Therefore, one can expect that the convolution of a binomial with a negative
binomial rv to be a more versatile discrete approximation, as it gives more flexibility in the choice of parameters to match the second moment, for example.

\begin{proposition}\label{binbp} Let $Y_{23} = Y_2 + Y_3$ have PGF  $ G_{23}(z)=G_2(z)G_3(z)$ and
 $p<q$. Then, for  $j\in\Zplus$ and  $g \in {\cal G}_{Y_{23}}$, the rv $Y_{23}$ has the following Stein operators:
\begin{eqnarray}
(\asta_1 g)(j)&=&(Mp+rq\qbar-pj+q\qbar j)g(j+1)+(r\qbar
p-Mp\qbar+p\qbar j)g(j+2)-qjg(j), \label{binb1}\\
(\asta_2
g)(j)&=&p\bigg(\frac{r\qbar}{p\pbar}+M-j\bigg)g(j+1)-qjg(j)\nonumber\\
&&+r(q\qbar+p)\sum_{m=2}^\infty\qbar^{m-1}\sum_{l=1}^{m-1}\Delta
g(j+l), \label{binb2}\\
(\asta_3
g)(j)&=&\qbar\bigg(\frac{Mp\pbar}{\qbar}+r+j\bigg)g(j+1)-jg(j)\nonumber\\
&&+M\bigg(\frac{p}{q}+\qbar\bigg)\sum_{m=2}^\infty(-1)^{m+1}
\bigg(\frac{p}{q}\bigg)^{m-1}\sum_{l=1}^{m-1}\Delta
g(j+l),\label{binb3}\\
(\asta_4
g)(j)&=&\bigg(Mp+\frac{r\qbar}{\pbar}\bigg)g(j+1)-jg(j)\nonumber\\
&&+\sum_{m=2}^\infty\bigg(M(-1)^{m+1}\bigg(\frac{p}{q}\bigg)^m+r\qbar^m\bigg)
\sum_{l=1}^{m-1}\Delta g(j+l). \label{binb4}
\end{eqnarray}\end{proposition}

\noindent \textbf{Proof. }
 Differentiating $G_{23}(z)= G_{2}(z)G_{3}(z) $ with respect to $z$, we obtain
\[
\sum_{k=0}^\infty\mu_kz^k\Big(\frac{Mp}{q+pz}+\frac{r\qbar}{1-\qbar
z}\Big)=\sum_{k=0}^\infty k\mu_kz^{k-1}.
\]
Multiplying both sides by $(q+pz)(1-\qbar z)$ and collecting the terms corresponding
to $z^k$, we obtain the recursive relation
\[\mu_k(Mp+r q\qbar-pk+q\qbar
k)+\mu_{k-1}(rp\qbar-Mp\qbar+p\qbar(k-1))-q\mu_{k+1}(k+1)=0.
\]
Multiplying the last equation by $g(k+1)$ and summing over all
nonnegative $k$, we obtain (\ref{binb1}). Observe next that
\[
\bigg(\frac{\pbar}{1-\qbar
z}\bigg)^r=\Exponent{r\sum_{i=1}^\infty\frac{\qbar^i}{i}(z^i-1)}.
\]
Therefore, (\ref{binb2})  follow from (\ref{a1}). Similarly,
(\ref{binb3}) follows from (\ref{a2}) and (\ref{uu1}), and
(\ref{binb4}) follows from (\ref{opCP}) and (\ref{uu1}).\qed

\begin{remark}
As discussed earlier, the operators $\asta_2$, $\asta_3$, and $\asta_4$ are binomial, negative binomial and Poisson perturbations, respectively. Note, however,  $\asta_1$ can not be seen as a perturbation operator.
\end{remark}

\subsection{Compound Distributions}
\noindent Next, we extend the $PGF$ technique for finding
Stein operators for a general class of compound distributions. 
Let
$S_N=\sum_{j=1}^N X_j$, where $N$ is a $\Zplus$-valued rv with $\mu_{k}=P(N=k)$
and the $X_j$ are iid rvs, independent of $N$, with $P(X_j=k)=p_k$ for $k \in \Zplus$. Here and henceforth, $S_0$ is treated as a degenerate rv
concentrated at zero. Then the $PGF$ of $S_N$ is
given by
\begin{equation*}
G_{S_N}(z)=G_N(G_{X_1}(z))=\sum_{j=0}^\infty\pi_jz^j,
\end{equation*}
where
\begin{equation}
\pi_j=P(S_N=j)=\sum_{k=0}^\infty
P(N=k)P(S_k=j)=\sum_{k=0}^\infty\mu_kp_{k,j}, \label{nine}
\end{equation}
and $p_{k,j}=P(S_k=j)$ denotes the $k$-fold convolution of
$\{p_j\}_{j\geqslant 0}$.  Thus,
\[
G_N(G_{X_1}(z))=\sum_{j=0}^\infty\bigg(\sum_{k=0}^\infty\mu_kp_{k,j}\bigg)z^j.\]
Further on, we assume that $\Expect (S_N)<\infty$. Then
\begin{equation}
G_{S_N}'(z)=\sum_{j=1}^\infty
j\pi_jz^{j-1}=\sum_{j=0}^\infty
(j+1)\pi_{j+1}z^j=\sum_{j=0}^\infty
(j+1)\bigg(\sum_{k=0}^\infty\mu_kp_{k,j+1}\bigg)z^j.
\label{twelve}
\end{equation}
Similarly,
\begin{eqnarray}
G_{S_N}'(z)
&=&\frac{d}{dG_{X_1}(z)}\sum_{k=0}^\infty\mu_k(G_{X_1}(z))^k\bigg(\frac{d}{dz}\sum_{m=0}^\infty
p_{1,m}z^m\bigg)\nonumber\\
&=&\sum_{k=0}^\infty (k+1)\mu_{k+1}(G_X(z))^k\sum_{m=0}^\infty
(m+1)p_{m+1}z^m. \label{fifteen}
\end{eqnarray}
Noting that $(G_{X_1}(z))^k=\sum_{s=0}^\infty p_{k,s}z^s$, we get
\begin{eqnarray}
G_{S_N}'(z)&=&\sum_{k=0}^\infty
(k+1)\mu_{k+1}\sum_{s=0}^\infty p_{k,s}z^s\sum_{m=0}^\infty
(m+1)p_{m+1}z^m\nonumber\\
&=&\sum_{s=0}^\infty\bigg\{\sum_{k=0}^\infty
(k+1)\mu_{k+1}\sum_{m=0}^sp_{k,m}(s-m+1)p_{s-m+1}\bigg\}z^s.
\label{eighteen}
\end{eqnarray}
Comparing (\ref{eighteen}) with (\ref{twelve}), we obtain the
required recursion relation, for $s\in\Zplus$, as
 \begin{equation}
(s+1)\sum_{k=0}^\infty\mu_kp_{k,s+1}=\sum_{k=0}^\infty
(k+1)\mu_{k+1}\sum_{m=0}^sp_{k,m}(s-m+1)p_{s-m+1}. \label{creq}
 \end{equation}
Next we derive a Stein operator. So far, some $\mu_j$'s were allowed
to be equal to zero. Now we restrict ourselves to the case
$\mu_j>0$, $j=0,1,2,\dots,K$ ($K=\infty$ is also allowed) and
assume that $\mu_{K+1}=\mu_{K+2}=\dots=0$, when $K < \infty$. Multiplying
(\ref{creq}) by $g(s+1)$ and summing over $s \in \Zplus$, we obtain
\[
\sum_{s=0}^\infty s g(s)\sum_{k=0}^K\mu_k p_{k,s}=
\sum_{s=0}^\infty g(s+1)\sum_{k=0}^K
(k+1)\mu_{k+1}\sum_{m=0}^sp_{k,m}(s-m+1)p_{s-m+1},
\]
or equivalently
\begin{equation*}
\sum_{k=0}^K\mu_k\sum_{m=0}^\infty
p_{k,m}\bigg(a_k\sum_{s=m}^\infty
g(s+1)(s-m+1)p_{s-m+1}-mg(m)\bigg)=0, \label{stop1}
\end{equation*}
where $a_k=(k+1)\mu_{k+1}/\mu_k$.  Changing the order of summation
in the above equation and setting $l=s-m+1$, we  obtain
\begin{equation}
\sum_{m=0}^\infty\sum_{k=0}^K\mu_kp_{k,m}\bigg(a_k\sum_{l=1}^\infty
g(l+m)lp_l-mg(m)\bigg)=0. \label{stop2}
\end{equation}
Next, let us assume that $a_k$'s satisfy Panjer's recursion:
$a_k=a+bk$ (see Panjer and Wang (1995)). From (\ref{nine}) and (\ref{stop2}),
\begin{equation}
\sum_{m=0}^\infty\pi_m\Big(a\sum_{l=1}^\infty
g(l+m)lp_l-mg(m)\Big)+b\sum_{m=0}^\infty\sum_{k=0}^K
k\mu_kp_{k,m}\sum_{l=1}^\infty g(l+m)lp_l=0. \label{stop3}
\end{equation}
Let $X$ be an independent copy of $X_1$. Then $\Expect
g(S_k+X)X=\Expect g(S_k+X)X_i$, $(i=1,2,\dots,k)$. Therefore,
\[
\sum_{m=0}^\infty kp_{k,m}\sum_{l=1}^\infty g(l+m)lp_l= k\Expect
g(S_k+X)X=\sum_{i=1}^k\Expect g(S_k+X)X_i=\Expect S_kg(S_k+X)
\]
and
\begin{eqnarray*}
\lefteqn{\sum_{m=0}^\infty\sum_{k=0}^K
k\mu_kp_{k,m}\sum_{l=1}^\infty g(l+m)lp_l
=\sum_{k=0}^K\mu_k\Expect
S_kg(S_k+X)}\hspace{2cm}\nonumber\\
&=& \sum_{k=0}^K\mu_k\sum_{m=0}^\infty mp_{k,m}\sum_{l=0}^\infty
g(l+m)p_l=\sum_{m=0}^\infty \pi_m m\sum_{l=0}^\infty g(l+m)p_l.
\end{eqnarray*}
Substituting the last expression into (\ref{stop3}), we  obtain a Stein operator as
\begin{equation}
(\asta g)(j)=\sum_{l=1}^\infty(al+bj)g(l+j)p_l-(1-bp_0)jg(j),\quad
j\in\Zplus.
 \label{stopmain}
\end{equation}
Thus, we have proved the following result.


\begin{theorem}\label{one} Let $N$ be rv concentrated on $\{0,1,2\dots,K\}$
($K$ may be infinite) with distribution $\mu_k=P(N=k)$
satisfying Panjer's recursion, for some $a,b \in {\mathbb R}$,
\[
\frac{(k+1)\mu_{k+1}}{\mu_k}=a+bk,\quad k=0,1,\dots,K,
\]
with $\mu_{K+1}=0$. Let $S_N=\sum_{j=1}^N X_j$, where the $X_j$ are
iid rvs independent of $N$ and concentrated on $\Zplus$ with
probabilities $P(X_1=k)=p_k$. If $\Expect (S_N)<\infty$ and $g \in {\cal G}_{S_N}$, then a
Stein operator for $S_N$ is given by (\ref{stopmain}).
\end{theorem}

\noindent \subsection {Some Examples}
\begin{itemize}
\item[a)]  Let  $N\sim P(\lambda)$, $\lambda >0$. Applying
Theorem \ref{one} with $K=\infty$, $a=\lambda$ and  $b=0$, we obtain
 \[
 (\asta g)(j)=\lambda\sum_{l=1}^\infty l
g(l+j)p_j-jg(j),\]
which coincides with the
one given in (\ref{opCP}) with $\lambda_j=\lambda p_j$.


\item[b)] Let $N \sim NB(r,\pbar)$, the negative binomial
distribution, $r>0$ and $0<\pbar<1$. Then $K=\infty$, $a=r\qbar$,
$b=\qbar$ and a Stein operator for the compound negative binomial distribution is
\begin{eqnarray}
(\asta
g)(j)&=&\qbar\sum_{m=1}^\infty(rm+j)g(j+m)p_m-(1-\qbar p_0)jg(j)\nonumber\\
&=&\sum_{m=1}^\infty
p_m\big\{\qbar(rm+j)g(j+m)-jg(j)\big\}-\pbar p_0jg(j)\nonumber\\
&=&\qbar(r\Expect X_1+j)-jg(j)\nonumber\\
&&-p_0\qbar j\Delta g(j)+\qbar\sum_{m=2}^\infty
(rm+j)p_m\sum_{k=1}^{m-1}\Delta g(j+k). \label{opCNB}
\end{eqnarray}
Note that the $PGF$ of $S_N$ is
\[
G_{S_N}(z)=\bigg(\frac{\pbar}{1-\qbar
G_X(z)}\bigg)^r=\bigg(\frac{\pbar}{1-\qbar\sum_{j=0}^\infty
p_jz^j}\bigg)^r.
\]
\item[c)] Let $N\sim Bi(n,p)$,  the binomial distribution, where $n\in\NN$ (the set of natural numbers) and $0<p<1$.
Then $K=n$, $a=np/q$, $b=-p/q$ and a Stein operator for the compound
binomial distribution is given by
\[
(\asta g)(j)=(p/q)\sum_{m=1}^\infty
(nm-j)g(j+m)p_m-(1+(p/q)p_0)jg(j)
\]
which can be written, in a form  similar to (\ref{astbin}), as
\begin{eqnarray}
(\asta g)(j)&=&p\sum_{m=1}^\infty
(nm-j)g(j+m)p_m-(q+pp_0)jg(j)\nonumber\\
&=&p(n\Expect X_1-j)g(j+1)-qjg(j)\nonumber\\
&&+pp_0j\Delta g(j)+\sum_{m=2}^\infty
(nm-j)p_m\sum_{k=1}^{m-1}\Delta g(j+k). \label{opCB}
\end{eqnarray}
Also, in this case
\[G_{S_N}(z)=(1+p(G_X(z)-1))^n=\bigg(1+p\sum_{j=0}^\infty
p_j(z^j-1)\bigg)^n.\]


\end{itemize}
\begin{remark}
(i) If we take $p_1=1$ in the examples above,  we obtain the standard Stein
operators for Poisson, binomial and negative binomial
distributions,  as given by (\ref{astpois}), (\ref{astbin}) and
(\ref{astnb}), respectively.

(ii)  Sometimes the form of $PGF$ allows to establish recursive relations
without differentiation. For example, the $PGF$ for the compound geometric
distribution is of the form
\[\frac{p}{1-q\sum_{m=1}^\infty p_mz^m}=\sum_{k=0}^\infty\mu_kz^k.
\]
Multiplying both sides by $1-q\sum_{m=1}^\infty p_mz^m$ and
collecting factors corresponding to $z^k$, we obtain
\[(\asta g)(j)=q\sum_{m=1}^\infty p_mg(j+m)-g(j).\]
This operator coincides with the one from Daly (2010). Note in this example $p_0=0$.
\end{remark}
\section{Perturbed Solutions to the Stein
Equation}\label{Sperturbations}
In this section, we discuss some known facts 
and explore properties of exact and approximate solutions to the Stein equation. Assume that $Y$ and $Z$ are rvs concentrated on $\Zplus$, $f \in {\cal F}$ and $g \in {\cal G}_Y$. 
Henceforth, $\norm{f}=\sup_{k}\ab{f(k)}$. As mentioned in Section \ref{Intro}, the second step in Stein's method is solving the equation (\ref{solo}).
%
Suppose a Stein operator for $Y$ is given by
\begin{equation}(\asta g)(j)=\alpha_j
g(j+1)-\beta_jg(j),\label{bx}\end{equation}
 where $\beta_0=0$ and $\alpha_k-\alpha_{k-1}\leqslant
\beta_k-\beta_{k-1}$  ($k=1,2,\dots$). Then a solution  $g$ to
(\ref{solo}) satisfies
\begin{equation}
\ab{\Delta g(j)}\leqslant 2\norm{f}
\min\bigg\{\frac{1}{\alpha_j},\frac{1}{\beta_j}\bigg\}. \quad
j\in\Zplus, f \in {\cal F}. \label{brxdelta}
\end{equation}
Define $g_i$ as a solution to (\ref{solo}) for the choice $f(j) = I(j=i)$, where $I(A)$ denotes the indicator function of $A$.
Then, from (2.18) and Theorem 2.10 of Brown and Xia (2001), we have
\begin{equation}
\ab{\Delta
g(i)}=\Ab{\sum_{j=0}^\infty f(j)\Delta
g_j(i)}\leqslant\sup_{j\geqslant 0}f(j)\ab{\Delta g_i(i)}\leqslant \sup_{j\geqslant 0}f(j) \min\{\alpha_i^{-1},\beta_i^{-1}\},\label{aux3}
\end{equation}
for  nonnegative functions $f$. The proof of (\ref{brxdelta}) can now be completed by following steps similar to that of  Lemma 2.2 from Barbour (1987), by noting the fact Stein equations with
$f^{+}(j)(:=f(j)-\inf_k f(k)\geqslant 0)$ and $f(j)$ on the right hand side of (\ref{solo}) have the same solution. If $f$ is non-negative, then $f^{+}(j)$  is not needed and $2\norm{f}$  in (\ref{aux3}) can be replaced by $\norm{f}$. Therefore, if $f: \Zplus\to [0,1]$, then $2\norm{f}$ in (\ref{brxdelta}) should be replaced by 1.

Note that different choices of $f$ lead to different probabilistic metrics. In
this paper, we consider total variation norm which is twice the
total variation metric. That is,
\begin{eqnarray*}
\norm{\eL(Y)-\eL(Z)}_{TV}&=&\sum_{j=0}^\infty\ab{P(Y=j)-P(Z=j)}=\sup_{\norm{f}\leqslant
1}\ab{\Expect f(Y)-\Expect f(Z)}\\
&=&2\sup_{f\in\F_1}\ab{\Expect f(Y)-\Expect
f(Z)}=2\sup_{A}\ab{P(Y\in A)-P(Z\in A)},
\end{eqnarray*}
where $\F_1=\{f| f:\Zplus\to [0,1]\}$, and the supremum is taken
over all Borel sets in the last equality.

Let $g$ be the solution to (\ref{solo}) for Poisson or negative
binomial or pseudo-binomial rv with Stein operator given by
(\ref{astpois}) or (\ref{astnb}) or (\ref{astbin}), respectively. Then the corresponding bounds are given respectively as
\begin{equation}
\norm{\Delta g}\leqslant\frac{2\norm{f}}{\max(1,\lambda)},\quad
\norm{\Delta g}\leqslant\frac{2\norm{f}}{r\qbar},\quad
\norm{\Delta g}\leqslant \frac{2\norm{f}}{\lfloor\widetilde N\rfloor
pq}. \label{Deltas}
\end{equation}
The first two bounds follow directly from (\ref{brxdelta}).
Observe that for pseudo-binomial distribution, the assumptions of
(\ref{brxdelta}) are not always satisfied. The last bound of
(\ref{Deltas}) follows from Lemma 9.2.1 in Barbour {\it et al.}~(1992b),
and using similar arguments as above.

If a Stein operator has a form different from
(\ref{bx}), then solving (\ref{solo}) and checking properties similar to (\ref{brxdelta}) becomes rather tedious.
 Apart
from the solution for compound geometric distribution by Daly (2010),
some partial success has been achieved for compound Poisson
distribution by Barbour and Utev (1998). In such situations, one can try the
perturbation technique introduced in Barbour and Xia
(1999) and further developed in Barbour and \v Cekanavi\v cius
(2002) and Barbour {\it et al.~}(2007). Roughly, the basic idea of perturbation can
be summarized in the following way: good properties of the
solution of (\ref{solo}) can be carried over to solutions
of Stein operators in similar forms.

Next, we formulate a partial case of Lemma 2.3 and Theorem 2.4 from Barbour {\it et al.~}(2007) under following setup.

Let $\asta_0$ be a Stein operator for  rv $Y$ with
support $\{0,1,2\dots, K\}$  ($K=\infty$ is allowed) and $g_0$ be the solution of the Stein equation
$$({\cal A}_0 g_0)(j) = f(j) - {\mathbb E}f(Y),~~~ f \in {\cal F}, ~ g_0 \in {\cal G}_Y.$$
Also, let there exist $\omega_1, \gamma >0$ such that $\norm{\Delta g_0}\leqslant
\omega_1\norm{f}\min(1,\gamma^{-1})$. Let ${\cal A}$ denote a Stein operator for rv $Z$ and
 $U := {\cal A} - {\cal A}_0$ be the perturbed part of ${\cal A}$ with respect to ${\cal A}_0$.

The following lemma establishes, under certain conditions, an approximation result between any two rvs $W$ and $Z$, using the observation that
a Stein operator for rv $Z$ can be seen as perturbation of a Stein operator for rv $Y$.
\begin{lemma} \label{lemma1}   Let  $Z$ be a rv with a Stein operator
$\asta=\asta_0+U$ and  $W$ be another rv, both concentrated on $\Zplus$. Also, assume that, for $g\in {\cal G}_Y \cap {\cal G}_Z$, there exist $\omega_2, \varepsilon >0$ such that
\[\norm{Ug}\leqslant \omega_2\norm{\Delta g}, \quad\hbox{}\quad \ab{\Expect(\asta
g)(W)}\leqslant\varepsilon\norm{\Delta g},\]
and $\omega_1\omega_2<\gamma$. Then
\[
\norm{\eL(W)-\eL(Z)}_{TV}\leqslant\frac{\gamma}{\gamma-\omega_1\omega_2}\Big(\varepsilon
\omega_1\min(1,\gamma^{-1})+2 P(Z > K) + 2P(W > K)\Big).
\]

\end{lemma}

\noindent
Next, using the assumptions of Lemma \ref{lemma1} and (\ref{Deltas}), we evaluate the values of $\omega_1$, $\omega_2$ and $\gamma$ to the various Stein operators derived in Section \ref{Soperator}. Our observations are as follows:

\begin{itemize}
\item[({\bf O1})] If a Stein operator is given by (\ref{opCP}), then we
have the Poisson perturbation with $\omega_1=2$, $\gamma=\sum_{m=1}^\infty
m\lambda_m$,
\[\norm{Ug}\leqslant \norm{\Delta g} \sum_{m=2}^\infty m(m-1)\ab{\lambda_m} =\norm{\Delta g} \omega_2
\] and
 $\omega_1\omega_2<\gamma$, provided $\{\lambda_m\}_{m \geqslant 2}$ is sufficiently small. For a general description of the problem, see Barbour {\it et~al.} (1992a).
\item[({\bf O2})] For the Stein operator given by (\ref{astaB}), we
have the pseudo-binomial perturbation with $\omega_1=2/pq$,
$\gamma=\lfloor M+\alpha/p\rfloor$, $\omega_2=p\alpha$ and
$\omega_1\omega_2<\gamma$,  if $p$ is sufficiently small (see Theorem \ref{BINO-POIS})
\item[({\bf O3})] Consider the Stein operator given by (\ref{astaCP}). Then we
have Poisson perturbation with $\omega_1=2$, $\gamma=Mp+\alpha$,
$\omega_2=Mp^2/(q-p)^2$ and $\omega_1\omega_2<\gamma$, whenever $p$ is
sufficiently small (see Theorem \ref{POIS-BIN}).
\item[({\bf O4})] For the Stein operator given by (\ref{binb2}), we
have the pseudo-binomial perturbation with $\omega_1=2$, $\gamma=\lfloor
M+r\qbar/(p\pbar)\rfloor pq$ and $\omega_2 =\frac{r\qbar(q\qbar+p)}{\pbar^2}$.
The condition $\omega_1\omega_2<\gamma$ is satisfied if $p$ and
$\qbar$ are sufficiently small.
\item[({\bf O5})] If the Stein operator is given by (\ref{binb3}), then we
have the negative binomial perturbation with $\omega_1=2$,
$\gamma=Mp\pbar+r\qbar$, $\omega_2=Mpq(p/q+\qbar)(q-p)^{-2}$ and
$\omega_1\omega_2<\gamma$, provided $p$ and $\qbar$ are sufficiently
small.
\item[({\bf O6})] Finally, consider the Stein operator given by (\ref{binb4}). Then we
have the Poisson perturbation, $\omega_1=2$, $\gamma=Mp+r\qbar/\pbar$,
$\omega_2=Mp^2/(q-p)^2+r\qbar^2/\pbar^2$ and
$\omega_1\omega_2<\gamma$, whenever $p$ and $\qbar$ are sufficiently
small.
\end{itemize}
\begin{remark}
(i) Note that, for the Stein operator in (\ref{binb1}), perturbation approach is not applicable. Also, for compound
negative binomial or compound binomial distributions, the perturbation part of the operator contains $j$, which makes the perturbation technique inapplicable, as the upper bound for $\|Ug\|$ can not be established. Consequently, either a new version of perturbation technique with nonuniform bounds should
be developed or a different approach should be devised.

(ii) We also remark here that once a Stein operator is derived (as discussed in Section \ref{Soperator}), the properties of the associated exact solution
to the Stein equation must be derived and this can be quite difficult. The perturbation approach, as discussed in some examples above (see {\rm (O1)-(O6)}), can be useful to get an upper bound on approximate solution to the Stein equation.
\end{remark}

\section{Application to Sums of Indicator Variables}\label{AIndicator}

In this section, we exploit the different forms of Stein operator to obtain better bounds for the approximation
problems to sums of possibly dependent indicator rvs.
In particular, we consider Stein operators derived in (\ref{astaB}) and (\ref{astaCP}) along with the corresponding observations (O2) and (O3) and establish the approximation results to the sums of independent and dependent indicators.

\noindent
Consider the sum $W=\sum_{i=1}^n \II_i$ of possibly dependent
indicator variables and let $W^{(i)}=W-\II_i$,
$P(\II_i=1)=p_i=1-P(\II_i=0)=1-q_i$ $(i=1,2,\dots,n)$.
Assume also $\widetilde W^{(i)}$
satisfy $P(\widetilde W^{(i)}=k)= P(W^{(i)}=k\vert \II_i=1)$, for all $k$.
 We choose $Y_{12}=Y_1+Y_2$ as the
approximating variable, where
$Y_1\sim P(\alpha)$,  $Y_2\sim Bi(M,p)$ and are independent. Denote its distribution by
$\BIP$ whose $ PGF$ is
given in (\ref{a3}). Poisson, signed compound Poisson and translated Poisson,
 binomial and  negative binomial approximations have been applied to the sums of
independent and dependent  Bernoulli variables in numerous papers;
see, for example, Barbour {\it et al.~}(1992b), Soon (1996), Barbour and
Xia (1999), Roos (2000), R\"{o}llin (2005), Pek\"{o}z {\it et al.~}(2009),
Daly {\it et al.~}(2012) and Vellaisamy {\it et al.~}(2013). Unlike asymptotic expansions or a signed
compound Poisson measure, $\BIP$ is a distribution. This might be an
added advantage in practical applications.

\subsection {The Choice of Parameters} Note that the $\BIP$ is a three-parametric
distribution. We choose the parameters $p$, $M$ and $\alpha$ to
ensure the almost matching of the first three moments of $W$. Denoting
as before the integral part by $\lfloor\cdot\rfloor$, we define
 \begin{eqnarray}
M&:=&\Bigg\lfloor{\bigg(\sivn p_i^2\bigg)^3} {\bigg(\sivn
p_i^3\bigg)^{-2}}\Bigg\rfloor,
 \label{pert1}\\
 \delta &:=&\bigg(\sivn p_i^2\bigg)^3\bigg(\sivn p_i^3\bigg)^{-2}-M,
\qquad 0\leqslant \delta <1,\label{pert2}\\
p&:=&\bigg(\sivn p_i^3\bigg)\bigg(\sivn p_i^2\bigg)^{-1}
;\qquad\alpha:=\sivn p_i-M p.
 \label{pert3}\end{eqnarray}
Then the following relations hold:
 \begin{equation}
Mp^2=\sivn p_i^2-\delta p^2,\quad
Mp^3=\sivn p_i^3-\delta p^3.
 \label{pert4}\end{equation}
Observe also that
\begin{equation*}
\bigg(\sivn p_i^2\bigg)^2\leqslant\sivn p_i\sivn p_i^3.
\end{equation*}
 Therefore, for
$\alpha>0$, the $\BIP$ is not a signed measure, but a distribution.
Similar to Soon (1996), we choose parameters to match the three moments
for the sum of independent Bernoulli variables. Thus, only weak
dependence of rvs is assumed. Note that the additional information
about dependence of rvs can significantly alter the choice of
parameters, see, for example, Daly {\it et al.~}(2011) and Corollary \ref{11bcp}. Observe that
$\alpha$ and $Mp$ can be of the same order.
 Indeed,
let $n$ be even and $p_1=p_2=\dots=p_{n/2}=1/6$,
$p_{n/2+1}=\dots=p_n=1/12.$.  Then $Mp=O(n)=\alpha$.

\subsection{Poisson Perturbation} We start with Stein operator  given in (\ref{astaCP}).  Some additional
notations are needed. Henceforth, let $\dirac_1$ and
$\dirac$ denote the degenerate distributions concentrated at 1 and 0,
respectively. The convolution operator is denoted by $*$. Also, let
\begin{eqnarray} \label{neqn39}
d&\isdefinedas&\Norm{\distributionof{W}{*}(\dirac_1-\dirac)^{*2}}_{TV}=
\sum_{k=0}^n\ab{\Delta^2P(W=k)}, \\ \label{neqn40}
d_1&\isdefinedas&\max_i\Norm{\distributionof{W^{(i)}}{*}(\dirac_1-\dirac)^{*2}}_{TV}=
\max_i\sum_{k=0}^n\ab{\Delta^2P(W^{(i)}=k)},\\  \nonumber
\widehat\lambda&=&\sum_{i=1}^np_i,\quad\sigma^2=\sum_{i=1}^n p_iq_i,\quad \tau=\max_ip_iq_i, \\  \nonumber
\eta_1&\isdefinedas& \sivn p_i(1+2p_i+4p_i^2)\Expect\ab{\widetilde
W^{(i)}-W^{(i)}}, \\  \label{theta1}
\theta_1&\isdefinedas& \frac{Mp^2}{(1-2p)^2(Mp+\alpha)}=\frac{\sum_{i=1}^np_i^2-\delta
p^2}{(1-2p)^2\sum_{i=1}^np_i}.
\end{eqnarray}
Now, we have the following $\BIP$ approximation result for the sum of weakly dependent indicator rvs.

\begin{theorem}\label{POIS-BIN}
Let $\max(p,\theta_1)<1/2$. Then
\begin{eqnarray*}
\norm{\distributionof{W}-\BIP}_{TV}&\leqslant&
\frac{2}{(1-2\theta_{1})\widehat\lambda}\bigg\{ d_1 \sivn p_i^4
+\frac{dMp^4}{(1-2p)^2}+(1+2p)\delta p^2+\eta_1 \bigg\}.
\label{PBtv}
\end{eqnarray*}
\end{theorem}
If the indicator variables are dependent, then obtaining the bounds for $d$ and
$d_1$ is difficult; see Lemma \ref{le1} and Daly (2011) for some partial cases and the history of the
problem. On the other hand, if the rvs are
independent, then by the unimodality of $W$ (see  Xia (1997)), we
obtain
\begin{equation*}
P(W=k)\leqslant\frac{1}{2\sigma},\qquad
 \norm{\distributionof{W}*(\dirac_1-\dirac)}_{TV}\leqslant
\frac{1}{\sigma}.\label{PWK}
\end{equation*}
 Now let $S_1$ and $S_2$ be the sets of indices
such that
\begin{equation*}
S_1\cup S_2=\{1,2,\dots,n\},\quad \sum_{i\in S_1}p_iq_i\geqslant
\frac{\sigma^2}{2},\quad \sum_{i\in S_2}p_iq_i\geqslant
\frac{\sigma^2-\tau}{2}.
\end{equation*}
Then, by the properties of total variation,
\begin{equation} d\leqslant\Norm{{\cal L}\Big(\sum_{i\in
S_1}\II_i\Big)*(\dirac_1-\dirac)}_{TV} \Norm{{\cal
L}\Big(\sum_{i\in
S_2}\II_i\Big)*(\dirac_1-\dirac)}_{TV}\leqslant\frac{2}{\sigma\sqrt{\sigma^2-\tau}}.
\label{PWK1}
\end{equation}
Similarly, \begin{equation} d_1\leqslant
\frac{2}{\sqrt{(\sigma^2-\tau)(\sigma^2-3\tau)}}. \label{PWK2}
\end{equation}
Thus, we have the following corollary for independent rvs.

\begin{corollary}\label{C-POIS-BIN}
Let $W$ be the sum of $n$ independent Bernoulli rvs with successs probabilities $p_i$, $\max(p,\theta_1)<1/2$ and
$\sigma^2>3\tau$. Then
\begin{eqnarray}
\lefteqn{\norm{\distributionof{W}-\BIP}_{TV}}\hspace{1cm}\nonumber\\
&\leqslant&
\frac{2}{(1-2\theta_{1})\widehat\lambda}\bigg\{\frac{2\sivn
p_i^4}{\sqrt{(\sigma^2-\tau)(\sigma^2-3\tau)}}
+\frac{2Mp^4}{(1-2p)^2\sigma\sqrt{\sigma^2-\tau}}+(1+2p)\delta p^2
\bigg\}. \label{CPBtv}
\end{eqnarray}
\end{corollary}

\begin{remark} (i) Observe that $\theta_1<p(1-2p)^{-2}\leqslant \max_i p_i(1-2\max_i p_i)^{-2}$.
Therefore, a
sufficient condition for $\max(p,\theta_1)<1/2$  is
$\max_ip_i<(3-\sqrt{5})/4=0.19098\dots ~.$

(ii) If all $p_i\asymp C$, then the order of accuracy of the bound in
(\ref{CPBtv}) is $O(n^{-1})$. In comparison to the Edgeworth
expansion,  the $\BIP$ is more advantageous since the approximation holds
for the total variation norm and no additional measures compensating
for the difference in supports are needed.

(iii) Also, one can compare (\ref{CPBtv}) with the classical Poisson approximation result (see, Chen and R\"{o}llin (2013) eq. (1.1)-(1.2)), where for $p_i \asymp C$ and the order of accuracy is $O(1)$.
\end{remark}

\noindent \Proofof{Theorem \ref{POIS-BIN}}Applying Newton's expansion,
similar to Barbour and \v Cekanavi\v cius (2002, p. 518),
  we get
\begin{eqnarray}
\lefteqn{\Ab{\Expect\Delta g(W+k)-\Expect\Delta
g(W+1)-(k-1)\Expect\Delta^2
g(W+1)}}\hspace{ 0.2cm}\nonumber\\
&&\leqslant\sum_{s=1}^{k-2}(k-1-s)\Ab{\Expect\Delta^3g(W+s)}\leqslant
\sum_{s=1}^{k-2}(k-1-s)\Ab{\sum_{j=0}^\infty\Delta
g(j+s)\Delta^2P(W=j-2)}\nonumber\\
&&\leqslant\frac{(k-1)(k-2)}{2}\norm{\Delta g}d. \label{aux1}
\end{eqnarray}
 By the definition of $M$ and $p$, defined respectively in (\ref{pert1}) and (\ref{pert3}),
 \begin{eqnarray}
 &&-M\sum_{l=2}^\infty\bigg(\frac{-p}{q}\bigg)^l(l-1)=-\sum_{k=1}^np_k^2+\delta
 p^2,\quad
 -M\sum_{l=2}^\infty\bigg(\frac{-p}{q}\bigg)^l\sum_{k=1}^{l-1}(k-1)=\sum_{k=1}^np_k^3-\delta
 p^3,\nonumber\\
&&~~~~~~~~M\sum_{k=2}^\infty\bigg(\frac{p}{q}\bigg)^l\sum_{k=1}^{l-1}(k-1)(k-2)=\frac{2Mp^4}{(1-2p)^4}.\label{aux2}
\end{eqnarray}
 Therefore, from (\ref{aux1}) and (\ref{aux2}), we get
\begin{eqnarray}
\lefteqn{\Ab{-M\sum_{l=2}^\infty\bigg(\frac{-p}{q}\bigg)^l
\sum_{k=1}^{l-1}\Expect\Delta g(W+k) +\sivn p_i^2\Expect\Delta
g(W+1)-\sivn p_i^3\Expect\Delta^2
g(W+1)}}\hspace{3cm}\nonumber\\
&\leqslant& \frac{Mp^4}{(1-2p)^4}\norm{\Delta g}d+\ab{\delta
p^2\Expect\Delta g(W+1)}+\ab{\delta
p^3\Expect\Delta^2 g(W+1)}\nonumber\\
&\leqslant& \frac{Mp^4}{(1-2p)^4}\norm{\Delta g}d+\delta
p^2(1+2p)\norm{\Delta g}. \label{pb1}
\end{eqnarray}
Taking into account (\ref{astaCP}) and (\ref{pb1}), we obtain
\begin{eqnarray}
\ab{\Expect(\asta g)(W)}&\leqslant&\Ab{\Expect\Big\{\sivn p_i
g(W+1)-Wg(W)\Big\}-\sum_{i=1}^np_i^2\Expect\Delta
g(W+1)\nonumber\\
&&+\sum_{i=1}^np_i^3\Expect\Delta^2g(W+1)}+\norm{\Delta
g}\bigg(\frac{Mp^4d}{(1-2p)^2}+\delta p^2(1+2p)\bigg)\nonumber
\\
&\leqslant& J_1+J_2+J_3+\norm{\Delta g}\bigg(\frac{Mp^4d}{(1-2p)^2}+\delta
p^2(1+2p)\bigg) (say). \label{pb11}
\end{eqnarray}
Here,
\begin{eqnarray}
J_1&=&\Ab{\Expect\Big\{\sivn p_i g(W+1)-Wg(W)\Big\}-\sivn
p_i^2\Expect\deltdui{1}{1}}\nonumber\\
&\leqslant&\Ab{\sivn
p_iq_i\Big(\Expect\left\{g(W^{(i)}+1)|\II_i = 0\right\} - \Expect\left\{g(W^{(i)}+1)|\II_i = 1\right\}\Big)}\nonumber\\
&=&\Ab{\sivn p_i\Expect\Big( g(W^{(i)}+1)- g(\widetilde
W^{(i)}+1)\Big)}\nonumber\\
&\leqslant& \norm{\Delta g}\sivn p_i\Expect\ab{W^{(i)}-\widetilde
W^{(i)}}.\label{pb2}
\end{eqnarray}
 Similarly,
\begin{eqnarray}
J_2&=&\Big\vert \sivn p_i^2\Expect\deltdui{1}{1}-\sivn
p_i^2\Expect\Delta g(W+1)\nonumber \\
&&+\sivn p_i^3
\Expect\cdik{2}{W^{(i)}+1}{1}\Big\vert\nonumber\\
&=&\Ab{\sivn
p_i^2q_i\Big(\Expect\deltdui{1}{0}-\Expect\deltdui{1}{1}\Big)}\nonumber\\
&\leqslant& 2\norm{\Delta g} \sivn p_i^2\Expect\ab{W^{(i)}-\widetilde W^{(i)}}\label{pb3}
\end{eqnarray}
and
\begin{eqnarray}
J_3&=&\Ab{\sivn p_i^3\Expect\Delta^2g(W+1)-\sivn
p_i^3\Expect\cdik{2}{W^{(i)}+1}{1}}
\nonumber\\
&\leqslant& \sivn
p_i^3q_i\ab{\Expect\cdik{2}{W^{(i)}+1}{0}-\Expect\cdik{2}{W^{(i)}+1}1}
\nonumber\\
&&+\sivn p_i^4\ab{\Expect\Delta^3 g(\widetilde W^{(i)}+1)}
\leqslant\norm{\Delta^3 g}\sivn p_i^3\Expect\ab{W^{(i)}-\widetilde
W^{(i)}}+\sivn p_i^4\norm{\Delta g}d_1.\label{pb4}
\end{eqnarray}
Collecting the bounds in (\ref{pb1})-(\ref{pb4}), applying Lemma
\ref{lemma1} and (O3) with $T=\infty$, the  proof is completed.
$\square$

\subsection{Binomial Perturbation}
Here, we approximate $W$ using Stein operator  in (\ref{astaB}).
In addition to the notations used above, let
 \begin{eqnarray*}
 d_2&\isdefinedas&\max_{i,j}\Norm{\distributionof{W^{(ij)}}{*}(\dirac_1-\dirac)
}_{TV}= \max_{i,j}\sum_k\ab{\Delta P(W^{(ij)}=k)},\\
 \widehat
T&\isdefinedas&\lfloor
M+\alpha/p\rfloor,\quad\theta_2\isdefinedas\frac{\alpha}{q\widehat
T},\quad W^{(ij)}=W-\II_i-\II_j.
\end{eqnarray*}
Also, let the distribution of
  $\widetilde W^{(ij)}_i$
  satisfy
  $P(\widetilde W^{(ij)}_i=k)=P(W^{(ij)}=k\vert\II_i=1)$, for
  all $k$.

\begin{theorem} \label{BINO-POIS} Let $\theta_2<1/2$. Then
\begin{eqnarray}
\lefteqn{\norm{\distributionof{W}-\BIP}_{TV}}\hspace{0.5cm}\nonumber\\
&\leqslant&\frac{2}{pq\widehat T(1-2\theta_2)}\bigg\{d_2\bigg(\sivn
p_i^4-p\sivn p_i^3\bigg)+\delta
p^2+\sivn p_i(2+2\ab{p_i-p})\EE\ab{\widetilde W^{(i)}-W^{(i)}}\nonumber\\
&&+\bigg(2\skn p_k^2\bigg)^{-1}\sum_{i,j=1}^n
p_ip_j\ab{p_i-p_j}\big[d_2\ab{p_i-p_j}
\ab{Cov(\II_i,\II_j)}+4 p_ip_j\Expect\ab{\widetilde
W^{(ij)}_i-\widetilde W^{(ij)}}\big]\bigg\}\nonumber\\
&&+\frac{2}{1-2\theta_2}\Big(P(Y_1+Y_2>\widehat T)+P(W>\widehat T)\Big).
\label{BPtv}
\end{eqnarray}
\end{theorem}
When  the indicator rvs are independent, a bound for the term $d_2$, similar to the one  in (\ref{PWK2}) for $d_1$,
can be obtained.
This leads to the following corollary.

\begin{corollary} \label{I-BINO-POIS} Let $W$ be the sum of $n$ independent indicator rvs, $\theta_2<1/2$ and $\sigma^2>3\tau$. Then
\begin{eqnarray}
\norm{\distributionof{W}-\BIP}_{TV}&\leqslant&
\frac{2}{1-2\theta_2}\bigg\{\frac{4}{pq\widehat
T(\sigma^2-3\tau)}\bigg(\sivn p_i^4-p\sivn
p_i^3\bigg)\nonumber\\&&+ \delta p^2+P(W>\widehat T)+P(Y_1+Y_2>\widehat
T)\bigg\}. \label{CBPtv}
\end{eqnarray}
\end{corollary}
%
\begin{remark} (i) If the rvs are independent, then
\[P(W>\widehat T)+ P(Y_1+Y_2>\widehat T)\leqslant \exponent{-\widehat\lambda\psi(p)},\]
where $\psi(p)=(p)^{-1}(-\ln p-1)+1$. Indeed,
\begin{eqnarray*}
&&P(Y_1+Y_2\geqslant \widehat T+1)\leqslant\ee^{-x(\widehat
T+1)}\Expect\ee^{xY_1}\Expect\ee^{xY_2}\leqslant
\ee^{-x(\widehat T+1)}\exponent{\alpha(\ee^x-1)}(q+p\ee^x)^M\\
&&~~~~\leqslant
\exponent{-x\widehat\lambda/p+(Mp+\alpha)(\ee^x-1)}\leqslant
\exponent{-\widehat\lambda(x/p+1-\ee^x)}.
\end{eqnarray*}
Now it suffices to take $x=-\ln p$. Similarly, one can obtain a bound for $P(W>\widehat T)$.
Observe, that $\psi(p)>0$ for any $p<1$.

(ii) If $p_i=C$, then the bound in (\ref{CBPtv}) is at least of the
order $O(n^{-1})$. The corresponding bounds for the binomial approximation as given in Corollary 1.3 of Soon (1996)  are of order $O(n)$ and  the ones in Remarks 2 of Roos (2000)  are of order $O(n^{-1/2})$. Also, see Theorem 1 of Ehm (1991)  where the bound is of order $O(1)$.

(iii) If all the $p_i$ are equal, then both sides of (\ref{CBPtv})
are equal to zero, as is the case for the binomial approximation  (see Soon (1996)).

(iv) Comparing Theorem \ref{BINO-POIS} with Theorem \ref{POIS-BIN}, we
observe that both have similar accuracy with respect to $\widehat \lambda$.
On the other hand, Theorem \ref{BINO-POIS} reflects the closeness
of $p_i$ and, in this sense, is more accurate than Theorem
\ref{POIS-BIN}.

(v) The $\BIP$ approximation (matching the first three moments) provides bounds with better accuracy (see Theorems \ref{POIS-BIN} and \ref{BINO-POIS}) than the bounds obtained (matching the two moments) for the binomial approximation (see Soon (1996) and Roos (2000)).
\end{remark}


\noindent \Proofof{Theorem \ref{BINO-POIS}} Using (\ref{astaB}) and (\ref{pert1})-(\ref{pert3}), note that
\begin{equation}
(\asta g)(j)=\sum_{i=1}^np_ig(j+1)-jg(j)-pj\Delta
g(j)+p\alpha\Delta g(j+1).\label{astaB1}
\end{equation}

Therefore,
\begin{eqnarray}
\Ab{\EE(\steinoperator g)(W)}&=&\Ab{\sivn p_i\EE g(W+1)-\sivn\EE
\II_ig(W)- p\sivn p_i \EE\deltdui{1}{1}\nonumber\\
&&+\Big(\sivn p_i(p-p_i)
 +\delta p^2\Big)\EE\Delta
g(W+1)}\leqslant\delta p^2\ab{\EE\Delta g(W+1)}\nonumber\\
&&+\Ab{\sivn p_i(p_i-p)\EE\deltdui{1}{1}+\sivn p_i(p-p_i)\EE\Delta
g(W+1)}\nonumber\\
 &&+\sivn p_iq_i\ab{\EE\deltdui{1}{0}-\EE\deltdui{1}{1}}\nonumber\\
 &=&R_1+R_2+R_3 \quad\hbox{(say).}\label{na1}
\end{eqnarray}
It is easy to check that
\begin{eqnarray}
\ab{R_1}&\leqslant&\norm{\Delta g}\delta p^2,\label{na2}\\
\ab{R_3}&=&\sivn p_i\ab{\EE \Delta g(W^{(i)}+1)-\EE g(\Delta \widetilde
W^{(i)}+1)}\leqslant 2 \norm{\Delta g}\sivn p_i\ewai,\label{na3}\\
\ab{R_2}&\leqslant&\Ab{\sivn
p_i^2(p-p_i)\EE\cdik{2}{W^{(i)}+1}{1}}\nonumber\\
&&+\sivn
p_iq_i\ab{p_i-p}\ab{\EE\deltdui{1}{0}-\EE\deltdui{1}{1}}.\label{a4}
\end{eqnarray}
The second summand in (\ref{a4}) is less than or equal to
\begin{equation} \norm{\Delta^2 g}\sivn
p_i\ab{p_i-p}\ewai\leqslant2\norm{\Delta g}\sivn
p_i\ab{p_i-p}\ewai. \label{a5}
\end{equation}
Also, the first term in (\ref{a4}) is
 \begin{equation}
    \sivn p_i^2(p_i-p)\EE \cdik{2}{W^{(i)}+ 1}{1}=\skn\Delta^2g(k)\sivn
    p_i(p_i-p)P(\II_i=1,W=k).
\label{a6}
\end{equation}
Moreover,
\begin{eqnarray}
\sivn p_i(p_i-p)P(\II_i=1,W=k)&=&\Big(\skn
p_k^2\Big)^{-1}\sum_{i,j}p_ip_j^2(p_i-p_j)P(\II_i=1,W=k)\nonumber\\
&=&\Big(2\skn
p_k^2\Big)^{-1}\bigg\{\sum_{i,j}p_ip_j^2(p_i-p_j)P(\II_i=1,W=k)\nonumber\\
&&+ \sum_{i,j}p_jp_i^2(p_j-p_i)P(\II_j=1,W=k)\bigg\}.\label{a7}
\end{eqnarray}

\noindent Set
\begin{equation} P_{ij}(k)=P(W=k+1\vert \II_i=1,\II_j=1)-P(W=k\vert
\II_i=1,\II_j=1). \label{a8}\end{equation}
Then it can be seen
(see Soon (1996), p.~709)) that
\begin{equation*}
P(\II_i=1,W=k)=P(\II_j=1,W=k)+(p_i-p_j)P(W^{(ij)}=k-1)+ Cov
\big(\II_i-\II_j,\II\{W^{(ij)}=k-1\}\big)
\end{equation*}
and
\begin{eqnarray*}
\lefteqn{p_iP(W^{(ij)}=k-1)-P(\II_i=1,W=k)}\hspace{2cm}\\
&=& \Big(p_ip_j+Cov\big(\II_i,\II_j\big)\Big)P_{ij}(k)-
Cov\big(\II_i,\II\{W^{(ij)}=k-1\}\big).
\end{eqnarray*}
Therefore,
 \begin{eqnarray*}
\lefteqn{\sivn p_i(p_i-p)P(\II_i=1,W=k)}\hspace{0.5cm}\\
&=&\frac{1}{2\sum_k
p_k^2}\bigg\{\sum_{i,j}p_i^2p_j(p_i-p_j)^2P(W^{(ij)}=k-1)-
\sum_{i,j}p_ip_j(p_j-p_i)^2P(\II_i=1,W=k)\\
&&+\sum_{i,j}p_i^2p_j(p_i-p_j)
Cov\big(\II_i-\II_j,\II\{W^{(ij)}=k-1\}\big)\bigg\}\\
&=& \frac{1}{2\sum_k
p_k^2}\bigg\{\sum_{i,j}p_i^2p_j^2(p_i-p_j)^2P_{ij}(k) +
\sum_{i,j}p_ip_j(p_j-p_i)^2Cov(\II_i,\II_j)P_{ij}(k)\\
&&- \sum_{i,j}p_ip_j(p_i-p_j)^2
Cov\big(\II_i,\II\{W^{(ij)}=k-1\}\big)\\
&&+ \sum_{i,j}
p_i^2p_j(p_i-p_j)Cov\big(\II_i-\II_j,\II\{W^{(ij)}=k-1\}\big)\bigg\}\\
&=& \frac{1}{2\sum_k
p_k^2}\bigg\{\sum_{i,j}p_i^2p_j^2(p_i-p_j)^2P_{ij}(k)+
\sum_{i,j}p_ip_j(p_i-p_j)^2Cov(\II_i,\II_j)P_{ij}(k)\\
&&+
\sum_{i,j}p_ip_j^2(p_i-p_j)Cov\big(\II_i,\II\{W^{(ij)}=k-1\}\big)\bigg\}.
\end{eqnarray*}
Consequently,
\begin{eqnarray}
\lefteqn{\Ab{\sivn p_i^2(p_i-p)\EE\{\Delta^2 g(W^{(i)}+1)\vert
\II_i=1\}}}\hspace{0.7cm}\nonumber\\
&\leqslant& \frac{1}{2\sum_kp_k^2}\bigg\{ \Ab{\sum_{k=1}^n\Delta^2
g(k)\sum_{i,j}p_i^2p_j^2 (p_i-p_j)^2P_{ij}(k)}\nonumber\\
&&+ \Ab{\sum_{k=1}^n\Delta^2 g(k)\sum_{i,j}p_ip_j(p_i-p_j)^2
Cov(\II_i,\II_j)P_{ij}(k)}\nonumber\\
&&+ 2\Ab{\sum_{k=1}^n\Delta^2 g(k)\sum_{i,j}p_ip_j^2(p_i-p_j)
Cov\big(\II_i, \II\{W^{(ij)}=k-1\}\big)}
\bigg\}\nonumber\\
&=&R_4+R_5+R_6 \quad\hbox{(say).} \label{a9}\end{eqnarray}
 We next derive upper bounds for $R_4$, $R_5$ and $R_6$ separately. First,
\begin{eqnarray}
R_4&\leqslant&
\frac{1}{2\sum_kp_k^2}\sum_{i,j}p_i^2p_j^2(p_i-p_j)^2
\Ab{\sum_{k=1}^n\Delta^2 g(k)
P_{ij}(k)}\nonumber\\
&\leqslant&
\frac{1}{2\sum_kp_k^2}\sum_{i,j}p_i^2p_j^2(p_i-p_j)^2\norm{\Delta
g} \sum_{k=1}^n\ab{\Delta P_{ij}(k-1)}\nonumber\\ &\leqslant& d_2
\norm{\Delta g}\Big(\sum_{k=1}^np_k^4-p\sivn p_i^3\Big).
\label{a10}\end{eqnarray}
Secondly,
 \begin{equation} R_5\leqslant
\frac{d_2\norm{\Delta
g}}{2\sum_kp_k^2}\sum_{i,j}p_ip_j(p_i-p_j)^2\ab{Cov(\II_i,\II_j)}.
\label{a11}\end{equation}
 Finally,
\begin{eqnarray*}
\lefteqn{Cov\big(\II_i,\II\{W^{(ij)}=k-1\}\big) = \EE
\II_i\II\{W^{(ij)}=k-1\}-p_i P(W^{(ij)}=k-1)} \\
&&= p_iP(W^{(ij)}=k-1\vert \II_i=1)- p_iP(W^{(ij)}=k-1).
\end{eqnarray*}
 Consequently,
\begin{eqnarray*}
\lefteqn{\Ab{\sum_{k=1}^n\Delta^2
g(k)\sum_{i,j}p_ip_j^2(p_i-p_j)Cov\Big(\II_i,
\II\{W^{(ij)}=k-1\}\Big)}}\hspace{1cm}\\
&=& \Ab{\sum_{i,j}p_i^2p_j^2(p_i-p_j)\Big(\EE\Delta^2 g(\widetilde
W^{(ij)}+1)- \EE\Delta^2
g(W^{(ij)}+1)\Big)}\\
&\leqslant& 4\norm{\Delta g}\sum_{i,j}p_i^2p_j^2\ab{p_i-p_j}\EE
\ab{\widetilde W^{(ij)}-W^{(ij)}}
\end{eqnarray*}
 and
\begin{equation} R_6\leqslant 2 \frac{\norm{\Delta
g}}{\sum_kp_k^2}\sum_{i,j}p_i^2p_j^2 \ab{p_i-p_j}\EE\ab{\widetilde
W^{(ij)}-W^{(ij)}}. \label{a12}\end{equation}
 Collecting the
bounds in (\ref{na1})-(\ref{a12}), we get the required bound for the
Stein operator defined in (\ref{astaB1}).
 Applying Lemma \ref{lemma1} and (O2),
the proof is completed. \qed


\subsection{Application to (1,1)-Runs}
We consider here a dependent setup arising out of independent Bernoulli trials. Let $\{X_j\}$ be a sequence of independent $Be(p^*)$ variables and $a(p^*)= p^*(1-p^*).$
Define, for $j \ge 2$,
\begin{equation}
{\mathbb I}_j = X_j (1 - X_{j-1})\;\;\mbox{ and }\;\; W = \sum_{j=2}^n {\mathbb I}_j.\label{IjW}
\end{equation}
Then, it can be easily seen that
\begin{eqnarray}
{\mathbb E}(W) &=& \sum_{j = 2}^n P({\mathbb I}_j = 1) =(n - 1)(1-p^*) p^* = (n-1) a(p^*) \; \mbox{(say)}, \label{mean}\\
{\mathbb V}(W) &=& (n - 1) a(p^*) + (5 - 3 n)(a(p^*))^2.\label{var}\\
\mbox{and}~~~ {\mathbb E}(W - {\mathbb E}W)^3 &=& (n-1)a(p^*) + (15 -9n)a(p^*)^2 + 4 (5n-11)a(p^*)^3.
\end{eqnarray}
This leads to the following choice of parameters:
\begin{eqnarray}
M &:=& \bigg\lfloor\frac{(3n -5)^2}{10 n - 22}\bigg\rfloor,\\
\delta &:=& \frac{(3n -5)^2}{10 n - 22} - M,~~~0\le \delta < 1,\\
p &:=& \left(\frac{10n -22}{3n -5}\right)a(p^*);~~~~~\alpha := (n -1)a(p^*) - M p.
\end{eqnarray}
\noindent
Let us define
\begin{equation}
K_1 = \frac{288(1 - 3 a(p^*))}{a(p^*)}\;\;\mbox{and}\;\; K_2 =
\frac{4}{ a(p^*)\sqrt{\min\{1 - a(p^*), 1/2\}}}. \label{K1K2}
\end{equation}
To apply Theorem 4.1, we need the following lemma.

\begin{lemma} \label{le1}
Let $\{{\mathbb I}_j\}_{j \ge 2}$ and $W$ be as defined in (\ref{IjW}), $d$ and $d_1$ be respectively defined in (\ref{neqn39}) and (\ref{neqn40}).
Then, for $(n-2)a(p^*) \ge 8$,
\begin{eqnarray}
d &\le& \frac{K_1}{n-1} + \frac{K_2}{\sqrt{n-1}} := \gamma(n-1),\label{d}\\
d_1 &\le& \frac{K_1}{n-2}+\frac{K_2}{\sqrt{n-2}}, \label{d_1}
\end{eqnarray}
where $K_1$ and $K_2$ are as defined in (\ref{K1K2}).
\end{lemma}

An application of Theorem 4.1 leads to the following corollary.

\begin{corollary}\label{11bcp}
Let $W$ be as defined in (\ref{IjW}) and $\theta_1$ be as defined in (\ref{theta1}). Assume $\max(p, \theta_1) \le 1/2$ and $(n-2)p^*(1-p^*) \ge 8$.
Then
\begin{eqnarray}
\|{\cal L}(W) - \BIP\|_{TV}&\le&\frac{2}{(1-2\theta_1) \widehat{\lambda}} \left\{\left(n a(p^*)^4 + \frac{Mp^4}{(1-2 p)^2}\right)\left(\frac{K_1}{n-2}+\frac{K_2}{\sqrt{n-2}}\right)\right.\nonumber\\
&&\left.+ (1+2p)\delta p^2 +  (n-1) C_1\right\},\label{IjWB}
\end{eqnarray}
where $C_1 = 2 \max\{1, 2(1 -a(p^*))\}a(p^*)(1+2a(p^*) + 4a(p^*)^2)(1-a(p^*)(1 -a(p^*)))$.
\end{corollary}
\begin{remark}
The bound given in (\ref{IjWB}) is of order $O(1)$ and comparable to the existing bounds for Poisson approximation given in Theorem 2.1 of Vellaisamy (2004). Also, it is an improvement over the bound given in Theorem 2.1 of Godbole (1993) which is of order $O(n)$.
\end{remark}

\noindent
\Proofof{Lemma \ref{le1}}
Let $\rho_0=0$ and define the stopping times
$$\rho_j = \min \left \{l > \rho_{j-1}| {\mathbb I}_l = 1 \right \}.$$
From Huang and Tsai (1991), the $T_j = \rho_j - \rho_{j-1}$ are iid having $PGF$
\begin{eqnarray*}
{\mathbb E}(z^T)= \frac{a(p^*) z^2}{1 - z + a(p^*)z^2}.
\end{eqnarray*}
Hence, ${\mathbb E}(T) = 1/a(p^*)$ and ${\mathbb V}(T) = {1 - 3 a(p^*)}/(a^{2}(p^*)).$
Observe now that $\rho_j = \sum_{i=1}^j T_i$ is the waiting time for $j$-th occurrence of ${\mathbb I}_l$.
Then it follows that the average number of occurrences in a sequence
$\{{\mathbb I}_j\}_{2 \le j \le n}$ is $({n-1})/{{\mathbb E}(T)} = (n-1) a(p^*).$
Suppose now $k = \lfloor (n-1) a(p^*) \rfloor +1$. Then $\rho_k =
\sum_{j=1}^k T_j$ and by Proposition 4.6 of Barbour and Xia (1999), we get
\begin{eqnarray*}
\|{\cal L}(\rho_k)*(I_1 - I)\|_{TV} \le \frac{2}{\sqrt{k
a(p^*) \min\{u_1, 1/2\}}},
\end{eqnarray*}
where $u_1 = 1 - (1/2)\|{\cal L}(T)*(I_1 - I)\|_{TV}$. Now, it can be
easily seen that $\|{\cal L}(T)*(I_1 - I)\|_{TV} = 2 a(p^*)$ which
implies \begin{eqnarray*} \|{\cal L}(\rho_k)*(I_1 - I)\|_{TV} \le \frac{2}{\sqrt{k a(p^*) \min\{1 - a(p^*), 1/2\}}} \le
\frac{2}{\sqrt{(n-1) (a(p^*))^2 \min\{1 - a(p^*), 1/2\}}}.
\end{eqnarray*}
Define maximal coupling (see Barbour {\it et~al.} (1992b, p. 254))
\begin{eqnarray}
2 P(\rho_k \neq \rho_k') = \|{\cal L}(\rho_k)*(I_1 - I)\|_{TV} \le \frac{2}{\sqrt{(n-1) (a(p^*))^2 \min\{1 -a(p^*), 1/2\}}}. \label{rhok}
\end{eqnarray}
Let now $\rho_k' = \sum_{j=1}^k T_j'$ such that $T_j$'s are iid and
$\rho_j' = \rho_{j-1}' + T_j'$ with $\rho_0'=0.$ Define now
\begin{eqnarray*}
{\mathbb I}_i = \left\{ \begin{array}{lll} 0,& \rho_{j-1}' < i < \rho_j'; 1
\le j \le k\\
1, & \rho_j' = i; 1 \le j \le k\\
{\mathbb I}_i,& \rho_k' < i.
\end{array}\right.
\end{eqnarray*}
Then, for $\rho_k \le (n-1)$ and $\rho_k = \rho_k'+1$, we have $W =
W'+1$. Hence,
\begin{equation}
P(W'+1 \neq W) \le P(\rho_k > n-1) + P(\rho_k \neq \rho_k'+1).\label{W(W+1)}
\end{equation}
Using Chebyshev's inequality, we get
$$P(\rho_k > n-1) \le \frac{{\mathbb V}(\rho_k)}{(n-1 - {\mathbb E}(\rho_k))^2}.$$
As seen earlier,
$${\mathbb E}(\rho_k) = \frac{k}{a(p^*)}; \; {\mathbb V}(\rho_k) = k \frac{1 - 3 a(p^*)}{(a(p^*))^2}.$$
Assume now, without loss of generality, $(n-1) a(p^*) \ge 8$. Then
\begin{eqnarray}
P(\rho_k > n-1) &\le& \frac{k (1 -3 a(p^*))}{((n-1) a(p^*) - k)^2}\nonumber\\
&\le& \frac{1.125 (1 - 3a(p^*))}{(n-1) a(p^*)(0.125)^2}\nonumber\\
&=& \frac{72 (1 - 3 a(p^*))}{(n-1) a(p^*)} = K_1/(n-1) (say).\label{rhokn}
\end{eqnarray}
Hence, we obtain from (\ref{rhok}), (\ref{W(W+1)}) and (\ref{rhokn})
$$d \le 2\|{\cal L}(W)*(I_1-I)\|_{TV} \le \frac{K_1}{n-1} + \frac{K_2}{\sqrt{n-1}}.$$
This proves (\ref{d}).

\noindent
Using similar arguments and the fact that  $T_j$'s are iid, (\ref{d_1}) immediately follows. \\

\noindent
\Proofof {Corollary \ref{11bcp}}
The bounds for $d$ and $d_1$ in Theorem 4.1 are given by Lemma \ref{le1}.
Next, to compute ${\mathbb E}|\widetilde W^{(i)} - W^{(i)}|$, construct the following two-dimensional stochastic process $\{(Z_l^{i1},Z_l^{i0})\}_{l\ge i}$ with initial state $(Z_i^{i1}, Z_i^{i0}) =
(1,0)$, where ${\cal L}(Z_l^{ij}) = {\cal L}({\mathbb I}_l |{\mathbb I}_i = j) $, for $j
= 0,1$, having following marginal distributions.

(i) For $l \ge i+2$,
\begin{eqnarray*}
P((Z_l^{i1}, Z_l^{i0}) = (0,0))) &=& 1 - a(p^*)\\
P((Z_l^{i1}, Z_l^{i0}) = (0,1))) &=& 0\\
P((Z_l^{i1}, Z_l^{i0}) = (1,0))) &=& 0\\
P((Z_l^{i1}, Z_l^{i0}) = (1,1))) &=& a(p^*).
\end{eqnarray*}

(ii) For $i < l \le i+ 1$,
\begin{eqnarray*}
P((Z_l^{i1}, Z_l^{i0}) = (0,0))) &=& 1 - \frac{a(p^*)}{1 - a(p^*)}\\
P((Z_l^{i1}, Z_l^{i0}) = (0,1))) &=& \frac{a(p^*)}{1 - a(p^*)}\\
P((Z_l^{i1}, Z_l^{i0}) = (1,0))) &=& 0\\
P((Z_l^{i1}, Z_l^{i0}) = (1,1))) &=& 0.
\end{eqnarray*}

\noindent Also, the joint distributions satisfy

(i) For $l = i$
\begin{eqnarray*}
P((Z_{l+1}^{i1}, Z_{l+1}^{i0}) = (0,0),(Z_l^{i1}, Z_l^{i0}) =
(0,0)) &=& 1 - 2 \frac{a(p^*)}{1 - a(p^*)}\\
P((Z_{l+1}^{i1}, Z_{l+1}^{i0}) = (0,1),(Z_l^{i1}, Z_l^{i0}) =
(0,0)) &=&  \frac{a(p^*)}{1 - a(p^*)}\\
P((Z_{l+1}^{i1}, Z_{l+1}^{i0}) = (0,0),(Z_l^{i1}, Z_l^{i0}) = (0,1))
&=& \frac{a(p^*)}{1 - a(p^*)},
\end{eqnarray*}
and zero otherwise.

(ii) For $l = i+1$,
\begin{eqnarray*}
P((Z_{l+1}^{i1}, Z_{l+1}^{i0}) = (0,0),(Z_l^{i1}, Z_l^{i0}) =
(0,0)) &=& 1 - (2 - a(p^*)) \frac{a(p^*)}{1 - a(p^*)}\\
P((Z_{l+1}^{i1}, Z_{l+1}^{i0}) = (1,1),(Z_l^{i1}, Z_l^{i0}) =
(0,0)) &=& a(p^*)\\
P((Z_{l+1}^{i1}, Z_{l+1}^{i0}) = (0,0),(Z_l^{i1}, Z_l^{i0}) = (0,1))
&=& \frac{a(p^*)}{1 - a(p^*)},
\end{eqnarray*}
and zero otherwise.

(iii) For $l \ge i+2 $,
\begin{eqnarray*}
P((Z_{l+1}^{i1}, Z_{l+1}^{i0}) = (0,0),(Z_l^{i1}, Z_l^{i0}) =
(0,0)) &=& 1 - 2 a(p^*)\\
P((Z_{l+1}^{i1}, Z_{l+1}^{i0}) = (1,1),(Z_l^{i1}, Z_l^{i0}) =
(0,0)) &=& a(p^*)\\
P((Z_{l+1}^{i1}, Z_{l+1}^{i0}) = (0,0),(Z_l^{i1}, Z_l^{i0}) = (1,1))
&=& a(p^*),
\end{eqnarray*}
and zero otherwise.

\noindent Let us now define the random variables
\begin{eqnarray*}
\zeta &=&\min \{ k - i | Z_k^{i1} = Z_k^{i0} = 1\},\;\;
\mbox{for}\;\; k \ge i\\
\mbox{and}\;\;\widetilde\zeta &=& \min\{i - k | Z_k^{i1} = Z_k^{i0} =
1\}, \;\; \mbox{for}\;\; i \le k.
\end{eqnarray*}
Due to symmetry of the stochastic process about $i$, we have
suppressed the index $i$. The distribution of $\zeta$ is given by
\begin{eqnarray*}
P(\zeta = k) &=& \left\{
\begin{array}{lcl}
a(p^*),& \mbox{for} & 2 \le k \le 3 \\
a(p^*) \left(\frac{1 - 2 a(p^*)}{1 -  a(p^*)}\right)^{k - 3}, & \mbox{for} &k
\ge 4. \end{array} \right.
\end{eqnarray*}
Therefore,
\begin{eqnarray*}
{\mathbb E}(\zeta) = a(p^*) + \frac{1}{a(p^*)}.
\end{eqnarray*}
Also, due to symmetry, we have $\zeta \stackrel{\cal L}{=} \widetilde
\zeta $.

\noindent Define now
\begin{eqnarray*}
W_l^{(i)} = \sum_{j = 2}^{i - \widetilde \zeta} Z_j^{i1} = \sum_{j = 2}^{i
-
\widetilde \zeta} Z_j^{i0},~~~~W_r^{(i)} = \sum_{j=i+\zeta}^n Z_j^{i1} = \sum_{j=i+ \zeta}^n Z_j^{i0}~~~\mbox{and}~~~
\xi^{i1} &=& \sum_{j = (i - \widetilde \zeta +1) \vee 2}^{(i+ \zeta
-1)\wedge n} Z_j^{i1} - Z_i^{i1}.
\end{eqnarray*}
Thus,
\begin{eqnarray*}
\widetilde W^{(i)} &=& W_l^{(i)} + W_r^{(i)} + \xi^{i1}.
\end{eqnarray*}


\noindent
Let now
\begin{eqnarray*}
{\mathbb I}_j' = \left\{ \begin{array}{lll} Z_j^{i1}, &\mbox{ with probability
}& a(p^*) \\
Z_j^{i0}, & \mbox{with probability } & 1 - a(p^*), \end{array} \right.
\end{eqnarray*}
and ${\mathbb I}_j'' \stackrel{\cal L}{=} {\mathbb I}_j$, but ${\mathbb I}_j''$ is independent of
$\{(Z_j^{i1},Z_j^{i0})| j \in [i - \widetilde \zeta, i +\zeta] \}$. Then
\begin{eqnarray*}
Z_j := \left\{
\begin{array}{lll}
{\mathbb I}_j'', & \mbox{if}~ j \in [i - \widetilde \zeta, i +\zeta]\\
{\mathbb I}_j',& \mbox{if}~ j > i+ \zeta \;\; \mbox{or} \;\; j < i - \widetilde \zeta.
\end{array}
\right.
\end{eqnarray*}
Define
\begin{eqnarray*}
\xi^{i} = \sum_{j = (i - \widetilde \zeta +1) \vee m}^{(i+ \zeta
-1)\wedge n} Z_j - Z_i ~~;~~
{W^{(i)}}' = W_l^{(i)} + W_r^{(i)} + \xi^i
\end{eqnarray*}
so that $W^{(i)} \stackrel{\cal L}{=} {W^{(i)}}'$. Note that
\begin{eqnarray*}
{\mathbb E}(\xi^i) \le {\mathbb E}(\zeta + \widetilde \zeta -1 ) =
\frac{2}{a(p^*)} + 2 a(p^*) -1\\
{\mathbb E}(\xi^{i1}) \le {\mathbb E}(\zeta + \widetilde \zeta -2)
=\frac{2}{a(p^*)} + 2 a(p^*) - 2.
\end{eqnarray*}
Therefore,
\begin{eqnarray*}
{\mathbb E}|\widetilde W^{(i)} - W^{(i)}| &\le& {\mathbb E}|\xi^i - \xi^{i1}| \\
&\le& a(p^*) \max\{2(1 -a(p^*),1)\}{\mathbb E}(\zeta + \widetilde \zeta -2) \\
&=& 2 \max\{2(1 -a(p^*)),1)\}(1 - a(p^*)(1 - a(p^*))).
\end{eqnarray*}
Thus, the bound given in Theorem 4.1 becomes
\begin{eqnarray*}
\|{\cal L}(W) - \BIP\|_{TV}&\le&\frac{2}{(1-2\theta_1) \widehat{\lambda}} \left\{\left(n a(p^*)^4 + \frac{Mp^4}{(1-2p)^2}\right)\left(\frac{K_1}{n-2}+\frac{K_2}{\sqrt{n-2}}\right)\right.\\
&&\left.+ (1+2p)\delta p^2 +  (n-1)C_1\right\}.
\end{eqnarray*}
This proves the corollary.

\noindent {\bf Acknowledgments}. The authors are grateful to the reviewers for several suggestions, critical comments
and especially for pointing out some relevant references related to this work.


\singlespacing

\bibliographystyle{plain}

\end{document}